\documentclass[12pt, leqno]{amsart}

\usepackage{float}
\usepackage{amssymb,amsmath,amsthm,amscd}
\usepackage{mathrsfs}
\usepackage{enumerate}
\usepackage{color}
\usepackage[vcentermath,enableskew]{youngtab}
\usepackage[all]{xy}

\newcommand{\C}{{\mathbf C}}

\newcommand{\Z}{{\mathbf Z}}
\newcommand{\B}{{\mathbf{B}}}
\newcommand{\A}{{\mathbf A}}
\newcommand{\F}{{\mathbf F}}

\newcommand{\V}{{\mathbf V}}

\newcommand{\seteq}{\mathbin{:=}}

\theoremstyle{plain}
\newtheorem{lemma}{Lemma}[section]
\newtheorem{prop}[lemma]{Proposition}
\newtheorem{theorem}[lemma]{Theorem}
\newcommand{\Prop}{\begin{prop}}
\newcommand{\enprop}{\end{prop}}
\newcommand{\Lemma}{\begin{lemma}}
\newcommand{\enlemma}{\end{lemma}}
\newcommand{\Th}{\begin{theorem}}
\newcommand{\enth}{\end{theorem}}
\newtheorem{corollary}[lemma]{Corollary}
\newcommand{\Cor}{\begin{corollary}}
\newcommand{\encor}{\end{corollary}}
\newtheorem{definition}[lemma]{Definition}

\newcommand{\Def}{\begin{definition}}
\newcommand{\edf}{\end{definition}}

\theoremstyle{definition}

\newtheorem{example}[lemma]{Example}

\newcommand{\g}{{\mathfrak{g}}}
\newcommand{\Uq}{{U_q(\mathfrak{q}(n))}}
\newcommand{\qn}{{\mathfrak{q}(n)}}

\newcommand{\isoto}[1][]{\mathop{\xrightarrow[#1]%
{{\raisebox{-.6ex}[0ex][-.6ex]{$\mspace{2mu}\sim\mspace{2mu}$}}}}}
\newcommand{\eq}{\begin{eqnarray}}
\newcommand{\eneq}{\end{eqnarray}}

\newcommand{\eqn}{\begin{eqnarray*}}
\newcommand{\eneqn}{\end{eqnarray*}}

\newcommand{\on}{\operatorname}

\newcommand{\bni}{\be[{\rm(i)}]}
\newcommand{\bna}{\be[{\rm(a)}]}

\newcommand{\QED}{\end{proof}}

\newcommand{\cl}{\colon}

\newcommand{\ba}{\begin{array}}
\newcommand{\ea}{\end{array}}

\newcommand{\bi}{\begin{enumerate}[{\rm(i)}]}

\newcommand{\eqsub}{\begin{subequations}\begin{eqnarray}}
\newcommand{\eneqsub}{\end{eqnarray}\end{subequations}}

\newcommand{\ol}{\overline}

\newcommand{\nc}{\newcommand}
\nc{\la}{\lambda}
\nc{\lam}{\lambda}
\nc{\U}[1][\g]{U_q(#1)}
\nc{\te}{\tilde{e}}
\nc{\tei}{\tilde{e}_i}
\nc{\tf}{\tilde{f}}
\nc{\tfi}{\tilde{f}_i}
\nc{\tU}{\widetilde U_q(\g)}
\nc{\tE}{\tilde{E}}
\nc{\tF}{\tilde{F}}

\nc{\tkone}{\tilde{k}_{\ol{1}}}
\nc{\teone}{\tilde{e}_{\ol{1}}}
\nc{\tfone}{\tilde{f}_{\ol{1}}}

\nc{\teibar}{\tilde{e}_{\ol{i}}} \nc{\tfibar}{\tilde{f}_{\ol{i}}}
\nc{\tki}{\tilde{k}_{\ol {i}}}

\nc{\BZ}{{\mathbb{Z}}}
\nc{\al}{\alpha}
\nc{\qs}{{q}}
\nc{\lan}{\langle}
\nc{\ran}{\rangle}
\nc{\re}{{\mathrm{re}}}
\nc{\wt}{\operatorname{wt}}
\nc{\Uf}[1][\g]{U^-_q(#1)}
\nc{\Ue}{U^+_q(\g)}
\nc{\eps}{\varepsilon}
\nc{\vphi}{\varphi}
\nc{\sphi}{\varphi^*}
\nc{\seps}{\varepsilon^*}

\nc{\nn}{\nonumber}
\def\max{{\mathop{\mathrm{max}}}}
\nc{\vp}{\varpi}
\nc{\cls}{{\operatorname{cl}}}
\nc{\Wt}{{\operatorname{Wt}}}
\nc{\Us}{U'_q(\g)}
\nc{\La}{\Lambda}
\nc{\ro}{{\rm(}}
\nc{\rf}{{\rm)}}
\nc{\norm}{{\mathrm{norm}}}
\nc{\qbox}{\quad\mbox}
\nc{\braid}{{\mathfrak{B}}}
\nc{\Ad}{\operatorname{Ad}}
\nc{\Aut}{\operatorname{Aut}}
\nc{\dt}[1]{\tilde{\tilde #1}}
\nc{\Sn}{S^{{\mathrm{norm}}}}
\nc{\aff}{{\mathrm{aff}}}
\nc{\rk}{{\mathrm{rk}}}
\nc{\tQ}{\widetilde{Q}}
\nc{\tP}{\widetilde{P}}
\nc{\tW}{\widetilde{W}}
\nc{\Dyn}{\mathrm{Dyn}}
\nc{\tD}{\widetilde{\Delta}}
\nc{\height}{{\operatorname{ht}}}
\nc{\bl}{\bigl}
\nc{\br}{\bigr}
\nc{\Hecke}{\mathrm{H}}
\nc{\HA}{\Hecke^{\mathrm{A}}}
\nc{\HB}{\Hecke^{\mathrm{B}}}
\nc{\K}{\mathrm{K}}
\newcommand{\scbul}{{\,\raise1pt\hbox{$\scriptscriptstyle\bullet$}\,}}
\nc{\vac}{{\phi}}
\nc{\Bt}{\B_\theta(\g)}
\nc{\be}{\begin{enumerate}}
\nc{\ee}{\end{enumerate}}
\nc{\low}{{\mathrm{low}}}
\nc{\upper}{{\mathrm{up}}}
\nc{\Zodd}{\Z_{\mathrm{odd}}}
\nc{\Ft}[1][n]{\mathbb{P}\mathrm{ol}_{#1}}
\nc{\Ftf}[1][n]{\widetilde{\mathbb{P}\mathrm{ol}}_{#1}}
\nc{\KA}{\on{K}^{\mathrm{A}}}
\nc{\KB}{\on{K}^{\mathrm{B}}}
\nc{\Res}{\on{Res}}
\nc{\Fc}[1][{n,m}]{\mathbf{F}_{#1}}
\nc{\tphi}{\tilde{\varphi}}
\nc{\CO}{\mathscr{O}}

\title[Quantum Queer Superalgebra and Crystal Bases]
{Quantum Queer Superalgebra  \\ and Crystal Bases}

\author[D. Grantcharov, J. H. Jung, S.-J. Kang, M. Kashiwara, M. Kim]{Dimitar Grantcharov$^{1}$, Ji Hye Jung$^{2}$, Seok-Jin
Kang$^{3}$, \\ Masaki Kashiwara$^{4}$, Myungho Kim$^{5}$}

\address{Department of Mathematics \\
         University of Texas at Arlington \\ Arlington, TX 76021, USA}

         \email{grandim@uta.edu}

\address{Department of Mathematical Sciences
         and
         Research Institute of Mathematics \\
         Seoul National University \\ Seoul 151-747, Korea}

         \email{jhjung@math.snu.ac.kr}

\address{Department of Mathematical Sciences
         and
         Research Institute of Mathematics \\
         Seoul National University \\ Seoul 151-747, Korea}

         \email{sjkang@math.snu.ac.kr}

\address{Research Institute for Mathematical Sciences \\
          Kyoto University \\ Kyoto 606-8502, Japan \\
          \& Department of Mathematical Sciences
         and
         Research Institute of Mathematics \\
         Seoul National University \\ Seoul 151-747, Korea}

         \email{masaki@kurims.kyoto-u.ac.jp}

\address{Department of Mathematical Sciences
         and
         Research Institute of Mathematics \\
         Seoul National University \\ Seoul 151-747, Korea}
         \email{mkim@math.snu.ac.kr}

\thanks{$^{1}$This work was partially supported by NSA grant H98230-10-1-0207 and by Max Planck Institute for Mathematics, Bonn.}
\thanks{$^{2}$This work was partially supported by BK21 Mathematical Sciences Division and by NRF Grant \# 2010-0010753.}
\thanks{$^{3}$This work was partially supported by KRF Grant \# 2007-341-C00001 and by National Institute for Mathematical Sciences (2010 Thematic Program, TP1004).}
\thanks{$^{4}$This work was partially supported by Grant-in-Aid for Scientific Research (B) 23340005,
Japan Society for the Promotion of Science.}
\thanks{$^{5}$This work was partially supported by KRF Grant \# 2007-341-C00001 and
by NRF Grant \# 2010-0019516.}

\begin{document}

\maketitle
\begin{abstract}
In this paper, we develop the crystal basis theory for the quantum
queer superalgebra $\Uq$. We define the notion of crystal bases,
describe the tensor product rule, and present the existence and
uniqueness of crystal bases for $\Uq$-modules in the category
$\mathcal{O}_{int}^{\ge 0}$.
\end{abstract}

\section{Introduction}
The crystal bases are one of the most prominent discoveries of the
modern combinatorial representation theory. Immediately after its
first appearance  in 1990 in  \cite{Kas90}, the crystal basis
theory developed rapidly and attracted considerable mathematical
attention. Many important and deep results for symmetrizable
Kac-Moody algebras have been established in the last 20 years
following the original works \cite{Kas90, Kas91, Kas93}. In
particular, an explicit combinatorial realization of crystal bases
for classical Lie algebras was given in \cite{KN}.

In contrast to the case of Lie algebras, the crystal base theory
for Lie superalgebras is still in its beginning stage. A major
difficulty in the superalgebra case arises from the fact that  the
category of finite-dimensional representations is in general not
semisimple. Nevertheless, there is an interesting category of
finite-dimensional $U_q(\mathfrak{g})$-modules which is semisimple
for the two super-analogues of the general linear Lie algebra
$\mathfrak{gl} (n)$: $\mathfrak{g}  = \mathfrak{gl} (m | n)$ and
$\mathfrak{g} = \mathfrak{q} (n)$. This is the category ${\mathcal
O}_{int}^{\geq 0}$ of representations that appear as
subrepresentations of tensor powers ${\bf V}^{\otimes N}$ of the
natural representation ${\bf V}$ of $U_q(\mathfrak{g})$. The
semisimplicity of  ${\mathcal O}_{int}^{\geq 0}$ is verified  in
\cite{BKK} for $\mathfrak{g} = \mathfrak{gl} (m|n)$ and in
\cite{GJKK} for $\mathfrak{g} = \mathfrak{q}(n)$.

The crystal basis theory of ${\mathcal O}_{int}^{\geq 0}$ for the
general linear Lie superalgebra $\mathfrak{g}= \mathfrak{gl} (m |n)$
was developed in \cite{BKK}. In this case the irreducible objects in
${\mathcal O}_{int}^{\geq 0}$ are indexed by partitions having
so-called $(m,n)$-hook shapes. This combinatorial description
enables us to index the crystal basis of any irreducible object
$V(\lambda)$ in ${\mathcal O}_{int}^{\geq 0}$ with highest weight
$\lambda$ by the set $B(Y)$ of semistandard tableaux $Y$ of shape
$\lambda$. In addition to the existence of the crystal basis, the
decompositions of $V(\lambda) \otimes {\bf V}$ and $B(Y) \otimes
{\bf B}$, where ${\bf B}$ is the crystal basis for ${\bf V}$, have
been found in \cite{BKK}.

In this paper we focus on the second super-analogue of the general
linear Lie algebra: the queer Lie superalgebra $\qn$. It has been
known since its inception that the representation theory of  $\qn$
is more complicated compared to the other classical Lie
superalgebra theories.  A distinguished feature of $\qn$ is that
any Cartan subsuperalgebra has a nontrivial odd part. As a result,
the highest weight space of any highest weight $\qn$-module has a
structure of a Clifford module. In particular, every
$\mathfrak{gl} (n)$-component of a finite-dimensional $\qn$-module
appears with multiplicity larger than one (in fact, a power of
two). Important results related to the representation theory of
$\qn$ include  the $\qn$-analogue of the celebrated Schur-Weyl
duality discovered by Sergeev in 1984 \cite{Ser}, and character
formulae for all simple finite-dimensional representations found
by Penkov and Serganova in 1997 \cite{PS}. The foundations of the
highest weight representation theory of the quantum     queer
superalgebra $U_q (\qn)$ have been established in \cite{GJKK}. An
interesting observation in \cite{GJKK} is that the classical limit
of a simple highest weight $U_q (\qn)$-module is a simple highest
weight $U(\qn)$-module or a direct sum of two highest weight
$U(\qn)$-modules.

In view of the above remarks, it is clear that developing a
crystal basis theory for the category ${\mathcal O}_{int}^{\geq
0}$ of $U_q (\qn)$ is a challenging problem. The purpose of this
paper is to announce the results that lead to a solution of this
problem. Take the base field to be $\C((q))$. Our main theorem is
the existence and uniqueness of the crystal bases of $U_q (\qn)$-
modules in ${\mathcal O}_{int}^{\geq 0}$. The proofs will appear
in full detail in a forthcoming paper. To overcome the challenges
described above, we modify the notion of a crystal basis and
introduce the so-called {\it abstract $\qn$-crystal}. To do so we
first define {\it odd Kashiwara operators} $\tilde{e}_{\ol{1}}$,
$\tilde{f}_{\ol{1}}$, and $\tilde{k}_{\ol{1}}$, where
$\tilde{k}_{\ol{1}}$ corresponds to an odd element in the Cartan
subsuperalgebra of $\qn$. Then, a {\it crystal basis} for  a
$U_q(\qn)$-module $M$ in the category ${\mathcal O}_{int}^{\geq
0}$ is a triple $(L,B,(l_b)_{b \in B})$, where the crystal lattice
$L$ is a free $\C [[q]]$-submodule of $M$, $B$ is a finite
$\mathfrak{gl} (n)$-crystal,  $(l_b)_{b \in B}$ is a family of
vector spaces such that $L / qL = \bigoplus_{b \in B} l_b$, with a
set of compatibility conditions for the action of the Kashiwara
operators imposed in addition. The definition of a crystal basis
leads naturally to the notion of an abstract $\qn$-crystal an
example of which is the $\mathfrak{gl}(n)$-crystal $B$ in any
crystal basis $(L,B,(l_b)_{b \in B})$. The modified notion of a
crystal allows us to consider the multiple occurrence of
$\mathfrak{gl} (n)$-crystals corresponding to a highest weight
$U_q(\qn)$-module $M$ in ${\mathcal O}_{int}^{\geq 0}$  as a
single $\qn$-crystal. It is worth noting that $M$ is not
necessarily a simple module and that the $\qn$-crystal $B$ of $M$
depends only on the highest weight $\lambda$  of $M$, hence we may
write $B = B(\lambda)$. In order to find the highest weight vector
of $B(\lambda)$, we use the action of the Weyl group on
$B(\lambda)$ and define odd Kashiwara operators
$\tilde{e}_{\ol{i}}$ and $\tilde{f}_{\ol{i}}$ for $i = 2,...,n-1$.
Then the highest weight vector of $B(\lambda)$ is simply the
unique vector  annihilated by the $2n-2$ Kashiwara operators
$\tilde{e}_{i}$ and $\tilde{e}_{\ol{i}}$. In addition to the
existence and uniqueness of the crystal basis of $M$, we establish
an isomorphism ${\bf B} \otimes B(\lambda) \simeq
\bigsqcup_{\lambda + \varepsilon_j: {\rm strict}} B(\lambda +
\varepsilon_j)$ and explicitly describe the highest weight vectors
of ${\bf B} \otimes B(\lambda)$ in terms of the even Kashiwara
operators $\tilde{f}_{i}$ and the highest weight vector of
$B(\lambda)$. We conjecture that the highest weight vectors of
$B(\lambda) \otimes {\bf B}$ can be found in an analogous way with
the aid of the odd Kashiwara operators $\tilde{f}_{\ol{i}}$.

\section{The quantum queer superalgebra}\label{sec:qn}

For an indeterminate $q$, let $\F=\C((q))$ be the field of formal
Laurent series in $q$ and let $\A=\C[[q]]$ be the subring of $\F$
consisting of formal power series in $q$. Let $P^{\vee} = \Z k_1
\oplus \cdots \oplus \Z k_n$ be a free abelian group of rank $n$ and
let ${\mathfrak h} = \C \otimes_{\Z} P^{\vee}$. Define the linear
functionals $\epsilon_i \in \mathfrak{h}^*$ by $\epsilon_i(k_j) =
\delta_{ij}$ $(i,j=1, \ldots, n)$ and set $P= \Z \epsilon_1 \oplus
\cdots \oplus \Z \epsilon_n$. We denote by $\alpha_i = \epsilon_i -
\epsilon_{i+1}$ the {\em simple roots}.

\Def The {\em quantum queer superalgebra
$U_q(\mathfrak{q}(n))$} is the superalgebra over $\F$ with 1
generated by $e_i$, $f_i$, $e_{\ol i}$, $f_{\ol i}$ $(i=1, \ldots,
n-1)$, $q^{h}$ $(h\in P^\vee)$,
$k_{\ol j}$ $(j=1, \ldots, n)$ with the following
defining relations.
\begin{align}
\allowdisplaybreaks
\nonumber & q^{0}=1, \ \ q^{h_1} q^{h_2} = q^{h_1 + h_2} \ \ (h_1,
h_2 \in P^{\vee}), \\
\nonumber & q^h e_i q^{-h} = q^{\alpha_i(h)} e_i \ \ (h\in P^{\vee}), \displaybreak[1]\\
\nonumber & q^h f_i q^{-h} = q^{-\alpha_i(h)} f_i \ \ (h\in P^{\vee}), \displaybreak[1]\\
\nonumber & q^h k_{\ol j} = k_{\ol j} q^h, \displaybreak[1]\\
\nonumber & e_i f_j - f_j e_i = \delta_{ij} \dfrac{q^{k_i - k_{i+1}} - q^{-k_i
+ k_{i+1}}}{q-q^{-1}}, \displaybreak[1]\\
\nonumber & e_i e_j - e_j e_i = f_i f_j - f_j f_i = 0 \quad \text{if} \ |i-j|>1, \displaybreak[1]\\
\nonumber & e_i^2 e_j -(q+q^{-1}) e_i e_j e_i  + e_j e_i^2= 0  \quad \text{if} \ |i-j|=1,\displaybreak[1]\\
\nonumber & f_i^2 f_j - (q+q^{-1}) f_i f_j f_i + f_j f_i^2 = 0  \quad \text{if} \ |i-j|=1,\displaybreak[1]\\
%\end{aligned}
%\end{equation}
%\begin{equation}
%\begin{aligned}
\nonumber & k_{\ol i}^2 = \dfrac{q^{2k_i} - q^{-2k_i}}{q^2 - q^{-2}}, \\
\nonumber & k_{\ol i} k_{\ol j} + k_{\ol j} k_{\ol i} =0 \ \ (i \neq j),
\\
\nonumber & k_{\ol i} e_i - q e_i k_{\ol i} = e_{\ol i} q^{-k_i}, \\
          & k_{\ol i} f_i - q f_i k_{\ol i} = -f_{\ol i} q^{k_i}, \\
\nonumber & e_i f_{\ol j} - f_{\ol j} e_i = \delta_{ij} (k_{\ol i}
q^{-k_{i+1}} - k_{\ol{i+1}} q^{-k_i}), \\
\nonumber & e_{\ol i} f_j - f_j e_{\ol i} = \delta_{ij} (k_{\ol i}
q^{k_{i+1}} - k_{\ol{i+1}} q^{k_i}), \\
%&\nonumber e_{\overline{i}}f_{\overline{j}}+
%f_{\overline{j}}e_{\overline{i}} = \delta_{ij} \Big(
%\frac{q^{k_i+k_{i+1}} - q^{-k_i-k_{i+1}}}{q-q^{-1}} \\
 %&  \nonumber \hs{20ex} +(q-q^{-1})k_{\overline{i}}k_{\overline{i+1}} \Big), \\
\nonumber &e_i e_{\ol i} - e_{\ol i} e_i = f_i f_{\ol i} - f_{\ol i} f_i = 0, \\
%\nonumber &e_{\overline{i}}^2 = - \frac{q-q^{-1}}{q+q^{-1}}e_i^2,\;
%f_{\overline{i}}^2 = \frac{q-q^{-1}}{q+q^{-1}}f_i^2,\\
\nonumber &e_i e_{i+1} - q e_{i+1}e_i =
e_{\overline{i}}e_{\overline{i+1}}+ q
e_{\overline{i+1}}e_{\overline{i}}, \\
\nonumber&q f_{i+1}f_i - f_i f_{i+1} =
f_{\overline{i}}f_{\overline{i+1}}+ q f_{\overline{i+1}}f_{\overline{i}}, \\
\nonumber & e_i^2 e_{\overline{j}} - (q+q^{-1})e_i e_{\overline{j}}
e_i + e_{\overline{j}} e_i^2= 0 \quad \text{if} \ |i-j|=1, \\
\nonumber & f_i^2 f_{\overline{j}} - (q+q^{-1})f_i f_{\overline{j}}
f_i + f_{\overline{j}} f_i^2=0 \quad \text{if} \ |i-j|=1.
\end{align}
\edf

The generators $e_i$, $f_i$ $(i=1, \ldots, n-1)$, $q^{h}$ ($h\in P^\vee$)
are regarded as {\em even} and $e_{\ol i}$, $f_{\ol
i}$ $(i=1, \ldots, n-1)$, $k_{\ol j}$ $(j=1, \ldots, n)$ are {\em
odd}. From the defining relations, it is easy to see that the even
generators together with $k_{\ol 1}$
%$e_{\ol 1}$, $f_{\ol 1}$,
generate the whole algebra $U_q(\mathfrak{q}(n))$.

The superalgebra $U_q(\mathfrak{q}(n))$ is a Hopf superalgebra with
the comultiplication $\Delta\cl U_q(\mathfrak{q}(n)) \to
U_q(\mathfrak{q}(n)) \otimes U_q(\mathfrak{q}(n))$ defined by
\begin{equation}
\begin{aligned}
& \Delta(q^{h})  = q^{h} \otimes q^{h}\quad\text{for $h\in P^\vee$,} \\
& \Delta(e_i)  = e_i \otimes q^{-k_i + k_{i+1}} + 1 \otimes e_i, \\
& \Delta(f_i)  = f_i \otimes 1 + q^{k_i - k_{i+1}} \otimes f_i, \\
& \Delta(k_{\ol 1}) =k_{\ol 1}\otimes q^{k_1}+ q^{-k_1} \otimes k_{\ol 1}.
%& \Delta(e_{\ol 1}) = e_{\ol 1} \otimes q^{k_1 + k_2} + 1 \otimes
%e_{\ol 1} \\
%& \qquad \quad - (q-q^{-1}) e_1 \otimes q^{k_2} k_{\ol 1}, \\
%& \Delta(f_{\ol 1}) = f_{\ol 1} \otimes 1 + q^{-k_1 - k_2} \otimes
%f_{\ol 1} \\
%& \qquad \quad + (q-q^{-1}) q^{-k_2} k_{\ol 1} \otimes f_1.
\end{aligned}
\end{equation}

Let $U^{+}$ (resp.\ $U^{-}$) be the subalgebra of
$U_q(\mathfrak{q}(n))$ generated by $e_i$, $e_{\ol i}$ $(i=1,
\ldots, n-1)$ (resp.\ $f_i$, $f_{\ol i}$ ($i=1, \ldots, n-1$), and
let $U^{0}$ be the subalgebra
generated by $q^{h}$ ($h\in P^\vee$) and $k_{\ol j}$ $(j=1, \ldots, n)$.
In \cite{GJKK}, it was
shown that the algebra $U_q(\mathfrak{q}(n))$ has the {\em
triangular decomposition}:
\begin{equation}
U^{-} \otimes U^{0} \otimes U^{+}\isoto
U_q(\mathfrak{q}(n)).
\end{equation}

Hereafter, a $U_q(\mathfrak{q}(n))$-module is understood as a
$U_q(\mathfrak{q}(n))$-supermodule.
A $U_q(\mathfrak{q}(n))$-module $M$ is called a {\em weight module}
if $M$ has a weight space decomposition $M=\bigoplus_{\mu \in P}
M_{\mu}$, where
$$M_{\mu}\seteq  \{ m \in M ; q^h m = q^{\mu(h)} m \ \ \text{for all} \ h
\in P^{\vee} \}.$$ The set of weights of $M$ is defined to be
$$\wt(M) = \{\mu \in P ; M_{\mu} \neq 0 \}.$$

\Def A weight module $V$ is called a {\em highest weight module
with highest weight $\la \in P$} if $V$ is generated
by a finite-dimensional $U^0$-module ${\mathbf v}_{\la}$ satisfying the
following conditions:
\bna
\item $e_i v = e_{\ol i} v =0$ for all $v \in \mathbf{v}_{\la}$,
$i=1, \ldots, n-1$,
\item $q^h v = q^{\la(h)} v$ for all $v\in {\mathbf v}_{\la}$, $h \in
P^{\vee}$.
\ee  \edf

\noindent There is a unique irreducible highest weight module with
highest weight $\la \in P$ up to parity change. We denote it by
$V(\la)$.

Set
\begin{equation*}
\begin{aligned}
P^{\ge 0} = & \{ \la = \la_1 \epsilon_1 + \cdots + \la_n \epsilon_n
\in P\, ; \, \la_j \in \Z_{\ge 0} \ \ \text{for all} \ j=1, \ldots, n \}, \\
\La^{+} = & \{\la = \la_1 \epsilon_1 + \cdots + \la_n \epsilon_n \in
P^{\ge 0}\, ; \, \la_{i} \ge \la_{i+1} \ \text{and} \ \la_{i}=\la_{i+1} \
 \text{implies}  \\
 & \quad \qquad \qquad \qquad \qquad \qquad \la_{i} = \la_{i+1} = 0 \ \text{for all} \ i=1,
 \ldots,n-1\}.
\end{aligned}
\end{equation*}
Note that each element $\la \in \La^{+}$ corresponds to a {\em
strict partition} $\la = (\la_1 > \la_2 > \cdots > \la_r >0)$. Thus
we will call $\la \in \La^{+}$ a strict partition.

We define $\mathcal{O}_{int}^{\ge 0}$ to be the category of
finite-dimensional weight modules such that $\wt(M) \subset P^{\ge
0}$ and $k_{\ol i} |_{ M_{\mu}}=0$ for any $i \in \{1, \ldots, n
\}$ and $\mu \in P^{\ge 0}$ satisfying $\langle k_{\ol i}, \mu
\rangle=0$. The fundamental properties of the category
$\mathcal{O}_{int}^{\ge 0}$ are summarized in the following
proposition.

\vskip 2ex

\Prop [\cite{GJKK}] \hfill
\bna
\item Every $U_q(\mathfrak{q}(n))$-module in
$\mathcal{O}_{int}^{\ge 0}$ is completely reducible.
\item Every irreducible object in $\mathcal{O}_{int}^{\ge 0}$ has the
form $V(\la)$ for some $\la \in \La^{+}$.
\ee
\enprop

\vskip 1ex

\section{Crystal bases}

Let $M$ be a $U_q(\mathfrak{q}(n))$-module in
$\mathcal{O}_{int}^{\ge 0}$. For $i=1, 2, \ldots, n-1$, we define
the {\em even Kashiwara operators} on $M$ in the usual way. That is,
for a weight vector $u \in M_{\la}$, consider the {\em $i$-string
decomposition} of $u$:
$$u =\sum_{k\ge 0} f_i^{(k)} u_k,$$
where $e_i u_k =0$ for all $k \ge 0$, $f_i^{(k)} = f_i^{k} / [k]!$,
$[k]=\dfrac{q^k - q^{-k}}{q - q^{-1}}$, $[k]! = [k][k-1] \cdots
[2][1]$, and we define the even Kashiwara operators $\tei$,
$\tfi$ $(i=1, \ldots, n-1)$ by
\begin{equation}
\begin{aligned}
& \tei u = \sum_{k \ge 1} f_i^{(k-1)} u_k, \qquad
 \tfi u = \sum_{k \ge 0} f_i^{(k+1)} u_k.
\end{aligned}
\end{equation}
On the other hand, we define the {\em odd Kashiwara operators}
$\tilde{k}_{\ol {1}}$, $\tilde{e}_{\ol {1}}$, $\tilde{f}_{\ol
{1}}$ by
\begin{equation}
\begin{aligned}
\tkone & = q^{k_1-1}k_{\ol 1}, \\
\teone & = - (e_1 k_{\ol 1} - q k_{\ol 1} e_1) q^{k_1 -1}, \\
\tfone & = - (k_{\ol 1} f_1 - q f_1 k_{\ol 1}) q^{k_2-1}.
\end{aligned}
\end{equation}

Recall that an abstract $\mathfrak{gl}(n)$-crystal is a set $B$
together with the maps
$\tei, \tfi\cl B \to B \sqcup \{0\}$, $\vphi_i,
\eps_i \cl B \to \Z \sqcup \{-\infty\}$ $(i=1, \ldots, n-1)$, and $\wt\cl
B \to P$ satisfying the conditions given in \cite{Kas93}.
In this paper, we say that
an abstract $\mathfrak{gl}(n)$-crystal is a {\em $\mathfrak{gl}(n)$-crystal}
if it is realized as a crystal basis of a finite-dimensional
integrable $U_q(\mathfrak{gl}(n))$-module.
In particular,
we have $\eps_i(b)=\max\{n\in\Z_{\ge0}\,;\,\tei^nb\not=0\}$
and $\vphi_i(b)=\max\{n\in\Z_{\ge0}\,;\,\tfi^nb\not=0\}$ for any $b$
in a $\mathfrak{gl}(n)$-crystal $B$.

\Def Let $M= \bigoplus_{\mu \in P^{\ge 0}} M_{\mu}$ be a
$U_q(\mathfrak{q}(n))$-module in the category
$\mathcal{O}_{int}^{\ge 0}$. A {\em crystal basis} of $M$ is a
triple $(L, B, l_{B}=(l_{b})_{b\in B})$, where
\bna
\item $L$ is a free $\A$-submodule of $M$ such that

\bni
\item $\F \otimes_{\A} L \isoto M$,

\item $L = \bigoplus_{\mu \in P^{\ge 0}} L_{\mu}$, where $L_{\mu} = L
\cap M_{\mu}$,

\item  $L$ is stable under the Kashiwara operators $\tei$,
$\tfi$ $(i=1, \ldots, \mbox{n-1})$, $\tkone$, $\teone$, $\tfone$.
\end{enumerate}

\item $B$ is a $\mathfrak{gl}(n)$-crystal together with
the maps $\teone, \tfone \cl B \to B \sqcup \{0\}$ such that

\bni
\item $\wt(\teone b) = \wt(b) + \alpha_1$, $\wt(\tfone) = \wt(b) -
\alpha_1$,

\item for all $b, b' \in B$, $\tfone b = b'$ if and only if $b = \teone b'$.
\end{enumerate}

\item $l_{B}=(l_{b})_{b \in B}$ is a family of non-zero $\C$-vector spaces
such that

\bni
\item $l_{b} \subset (L/qL)_{\mu}$ for $b \in B_{\mu}$,

\item  $L/qL = \bigoplus_{b \in B} l_{b}$,

\item $\tkone l_{b} \subset l_{b}$,
\item for $i=1, \ldots, n-1, \ol 1$, we have
\be[{\rm(1)}]
\item
if $\tei b=0$ then $\tei l_{b} =0$, and
otherwise $\tei$ induces an isomorphism
$l_{b}\isoto l_{\tei b}$.
\item
if $\tfi b=0$ then $\tfi l_{b} =0$,
and otherwise $\tfi$ induces an isomorphism $l_{b}\isoto l_{\tfi b}$.
\ee
\end{enumerate}
\end{enumerate}

\edf

Note that one can prove that
 $\teone^2 = \tfone^2 = 0$ as endomorphisms of $L/qL$
for any  crystal basis $(L, B, l_{B})$.

\begin{example}
Let $$\V = \bigoplus_{j=1}^n \F v_{j} \oplus \bigoplus_{j=1}^n \F
v_{\ol j}$$ be the vector representation of $U_q(\mathfrak{q}(n))$.
The action of $\Uq$ on $\V$ is given as follows: \\
$e_iv_j=\delta_{j,i+1}v_i$, $e_iv_{\ol j}=\delta_{j,i+1}v_{\ol i}$,
$f_iv_j=\delta_{j,i}v_{i+1}$, $f_iv_{\ol j}=\delta_{j,i}v_{\ol{i+1}}$,
$e_{\ol i}v_j=\delta_{j,i+1}v_{\ol{i}}$,
$e_{\ol i}v_{\ol j}=\delta_{j,i+1}v_{i}$,
$f_{\ol i}v_j=\delta_{j,i}v_{\ol{i+1}}$,
$f_{\ol i}v_{\ol j}=\delta_{j,i}v_{{i+1}}$,
$q^h v_j=q^{\epsilon_j(h)} v_j$, $q^h v_{\ol j}=q^{\epsilon_j(h)} v_{\ol j}$,
$k_{\ol i}v_j=\delta_{j,i}v_{\ol j}$,
$k_{\ol i}v_{\ol j}=\delta_{j,i}v_{j}$. \\

Set $$\mathbf{L} = \bigoplus_{j=1}^n \A v_{j} \oplus
\bigoplus_{j=1}^n \A v_{\ol j},$$ $l_{j} = \C v_{j} \oplus \C
v_{\ol j}$, and let $\B$ be the crystal graph given below.

$$\xymatrix@C=5ex
{*+{\young(1)} \ar@<0.1ex>[r]^-{1}
\ar@{-->}@<-0.9ex>[r]_{\ol 1} & *+{\young(2)} \ar[r]^2 & *+{\young(3)} \ar[r]^3 & \cdots \ar[r]^{n-1} & *+{\young(n)} }$$

Here, the actions of $\tfi$ $(i=1, \ldots, n-1, \ol 1)$ are
expressed by $i$-arrows. Then $(\mathbf{L}, \B,
l_{\B}=(l_j)_{j=1}^n)$ is a crystal basis of $\V$.

\end{example}

\vskip 2ex

\Th Let $M_j$ be a $U_q(\g)$-module in $\mathcal{O}_{int}^{\ge
0}$ with crystal basis $(L_j, B_j, l_{B_j})$ $(j=1,2)$. Set $B_1
\otimes B_2 = B_1 \times B_2$ and
$$l_{B_{1} \otimes B_{2}}=(l_{b_1} \otimes l_{b_2})_{b_1 \in B_1,
b_2 \in B_2}.$$ Then $$(L_1 \otimes_{\A} L_2, B_1 \otimes B_2,
l_{B_1 \otimes B_2})$$ is a crystal basis of $M_1 \otimes_{\F} M_2$,
where the action of the Kashiwara operators on $B_1 \otimes B_2$ are
given as follows.

\begin{equation} \label{eq1:tensor product}
\begin{aligned}
\tei(b_1 \otimes b_2) & = \begin{cases} \tei b_1 \otimes b_2 \ &
\text{if} \ \vphi_i(b_1) \ge \eps_i(b_2), \\
b_1 \otimes \tei b_2 \ & \text{if} \ \vphi_i(b_1) < \eps_i(b_2),
\end{cases} \\
\tfi(b_1 \otimes b_2) & = \begin{cases} \tfi b_1 \otimes b_2 \
& \text{if} \  \vphi_i(b_1) > \eps_i(b_2), \\
b_1 \otimes \tfi b_2 \ & \text{if} \ \vphi_i(b_1) \le \eps_i(b_2),
\end{cases}
\end{aligned}
\end{equation}
\begin{equation} \label{eq2:tensor product}
\begin{aligned}
\teone (b_1 \otimes b_2) & = \begin{cases} \teone b_1 \otimes b_2
& \text{if $\lan k_1, \wt b_2 \ran = 0$,\ }  \lan k_2, \wt b_2 \ran =0, \\
b_1 \otimes \teone b_2   %& \text{if} \ \teone b_2 \neq 0, \\
%0 \ \
&  \text{otherwise,}
\end{cases} \\
\tfone(b_1 \otimes b_2) & = \begin{cases} \tfone b_1 \otimes b_2
& \text{if $\lan k_1, \wt b_2 \ran = 0$,\ } \lan k_2, \wt b_2 \ran =0, \\
b_1 \otimes \tfone b_2   %& \text{if} \ \tfone b_2 \neq 0, \\
%0 \
& \text{otherwise}.
\end{cases}\\
\end{aligned}
\end{equation}
 \enth
\begin{proof}[Sketch of Proof]
Our assertion follows from the following comultiplication formulas.
\begin{equation*}
\begin{aligned}
& \Delta(\tkone) = \tkone \otimes q^{2 k_1} + 1 \otimes \tkone, \\
& \Delta(\teone) = \teone \otimes q^{k_1 + k_2} + 1 \otimes \teone
 - (1-q^2) \tkone \otimes e_1 q^{2 k_1}, \\
& \Delta(\tfone) = \tfone \otimes q^{k_1 + k_2} + 1 \otimes \tfone
 - (1-q^{2}) \tkone \otimes f_1 q^{k_1 + k_2-1}.
\end{aligned}
\end{equation*}
\end{proof}

Motivated by the properties of crystal bases, we introduce the
notion of abstract crystals.

\Def An {\em abstract $\mathfrak{q}(n)$-crystal} is a
$\mathfrak{gl}(n)$-crystal together with the maps $\teone, \tfone\cl B
\to B \sqcup \{0\}$ satisfying the following conditions:
\bna
\item $\wt(B)\subset P^{\ge0}$,
\item $\wt(\teone b) = \wt(b) + \alpha_1$, $\wt(\tfone) = \wt(b) -
\alpha_1$,

\item for all $b, b' \in B$, $\tfone b = b'$ if and only if $b = \teone b'$.

\end{enumerate}
 \edf

Let $B_1$ and $B_2$ be abstract $\qn$-crystals. The {\em tensor
product} $B_1 \otimes B_2$ of $B_1$ and $B_2$ is defined to be the
$\mathfrak{gl}(n)$-crystal $B_1 \otimes B_2$ together with the maps
$\teone$, $\tfone$ defined by \eqref{eq2:tensor product}.
Then it is an abstract $\qn$-crystal. Note that $\otimes$ satisfies
the associative axiom.

\vskip 2ex

\begin{example}

\begin{enumerate}
\item If $(L, B, l_{B})$ is a crystal basis of a $\Uq$-module $M$ in the
category $\mathcal{O}_{int}^{\ge 0}$, then $B$ is an abstract
$\qn$-crystal.

\item The crystal graph $\B$ is an abstract $\qn$-crystal.

\item By the tensor product rule, $\B^{\otimes N}$ is an abstract
$\qn$-crystal. When $n=3$, the $\qn$-crystal structure of $\B
\otimes \B$ is given below.

$$\xymatrix
{*+{\young(1) \otimes \young(1)} \ar[r]^1 \ar@{-->}[d]^{\ol 1} &
 *+{\young(2) \otimes \young(1)} \ar@<-0.5ex>[d]_1 \ar@{-->}@<0.5ex>[d]^{\ol 1} \ar[r]^2&
 *+{\young(3) \otimes \young(1)} \ar@<-0.5ex>[d]_1 \ar@{-->}@<0.5ex>[d]^{\ol 1} \\
 *+{\young(1) \otimes \young(2)} \ar[d]^2 &
 *+{\young(2) \otimes \young(2)} \ar[r]_2 &
 *+{\young(3) \otimes \young(2)} \ar[d]^2 \\
 *+{\young(1) \otimes \young(3)} \ar@<-0.5ex>[r]_1 \ar@{-->}@<0.5ex>[r]^{\ol 1} &
 *+{\young(2) \otimes \young(3)} &
 *+{\young(3) \otimes \young(3)}
 }$$

\item For a strict partition $\la = (\la_1 > \la_2 > \cdots > \la_r
>0)$, let $Y_{\la}$ be the skew Young diagram having $\la_1$ many
boxes in the first diagonal, $\la_2$ many boxes in the second
diagonal, etc. For example, if $\la$ is given by $(7 > 6 > 4 > 2 > 0)$, then we have

$$Y_{\la} = \young(::::::\hfill,:::::\hfill\hfill,::::\hfill\hfill\hfill,:::\hfill\hfill\hfill\hfill,::\hfill\hfill\hfill\hfill,:\hfill\hfill\hfill,\hfill\hfill) \quad.$$

Let $\B(Y_{\la})$ be the set of all semistandard tableaux of shape
$Y_{\la}$ with entries from $1, 2, \ldots, n$. Then by an {\em
admissible reading} introduced in \cite{BKK}, $\B(Y_{\la})$ is
embedded in $\B^{\otimes |\la|}$ and it is stable under $\tei,\tfi$,
$\teone, \tfone$. Hence it becomes an abstract $\qn$-crystal.
Moreover, the $\qn$-crystal structure thus obtained does not depend
on the choice of admissible readings.
\end{enumerate}
\end{example}

\vskip 2ex

Let $B$ be an abstract $\qn$-crystal. For $i=2, \ldots, n-1$, let
$w$ be an element of the Weyl group $W$ with shortest length such
that $w(\alpha_i) = \alpha_1$. Such an element is unique and we
may choose $w=s_2 \cdots s_{i} s_1 \cdots s_{i-1}$. We define the
{\em odd Kashiwara operators} $\teibar$, $\tfibar$ $(i=2, \ldots,
n-1)$ by
$$\teibar = S_{w^{-1}} \teone
S_{w}, \ \ \tfibar = S_{w^{-1}} \tfone S_{w}.$$
Here $S_w$ is the Weyl group action on the $\mathfrak{gl}(n)$-crystal.
The operators
$\teibar$, $\tfibar$ do not depend on the choice of reduced
expressions of $w$. We say that $b\in B$ is a {\em highest weight
vector} if $\tei b = \teibar b =0$ for all $i=1, \ldots, n-1$.

\vskip 2ex

\section{Existence and uniqueness}

In this section, we present the main result of our paper.

\vskip 1ex

\Th
\bna
\item  Let $\la \in \La^{+}$ be a strict partition %with $|\la| \le N$
and let $M$ be a highest weight $\Uq$-module in the category
${\mathcal O}_{int}^{\geq 0}$ with highest weight $\la$. If $(L,
B, l_{B})$ is a crystal basis of $M$, then $L_\la$ is invariant
under $\tki\seteq q^{k_i-1} k_{\ol i}$ for all $i=1, \ldots, n$.
Conversely, if $M_{\la}$ is generated by a free $\A$-submodule
$L_{\la}^{0}$ invariant under $\tki$ $(i=1, \ldots, n)$, then
there exists a unique crystal basis $(L, B, l_{B})$ of $M$ such
that

\bni
\item $L_{\la}=L_{\la}^{0}$,

\item $B_{\la} = \{b_{\la} \}$,

\item $L_{\la}^{0} / q L_{\la}^{0} = l_{b_{\la}}$,

\item $B$ is connected.
\end{enumerate}

Moreover, as an abstract $\qn$-crystal, $B$
 depends only on $\la$. Hence we may write $B=B(\la)$.

\item The $\qn$-crystal $B(\la)$ has a unique highest weight vector
$b_{\la}$. % for $|\la| \le N$.

\item If $b \in \B \otimes B(\la)$ is a highest weight vector,
% with $|\la| \le N-1$,
then we have $$b = 1 \otimes \tf_{1} \cdots \tf_{j-1} b_{\la}$$
for some $j$ such that $\la + \epsilon_j$ is a strict partition.

\item
Let $M$ be a $U_q(\mathfrak{q}(n))$-module in
$\mathcal{O}_{int}^{\ge 0}$, and let $(L,B,l_B)$ be a crystal basis of $M$.
Then there exist decompositions
$M=\bigoplus_{a\in A}M_a$ as a $U_q(\mathfrak{q}(n))$-module,
$L=\bigoplus_{a\in A}L_a$ as an $\A$-module,
$B=\bigsqcup_{a\in A}B_a$ as a $\qn$-crystal, parametrized by a set $A$
such that for any $a\in A$ the following conditions hold:
\bni
\item
$M_a$ is a highest weight module with highest weight $\la_a$ and
$B_a \simeq B(\la_a)$ for some strict partition $\la_a$,
\item
$L_a=L\cap M_a$, $L_a/qL_a=\bigoplus_{b\in B_a}l_b$,
\item
$(L_a,B_a,l_{B_a})$ is a crystal basis of $M_a$.
\ee
\item Let $M$ be a highest weight $\Uq$-module
in the category ${\mathcal O}_{int}^{\geq 0}$ with highest weight
$\la$. % with $|\la| \le N-1$.
Assume that $M$ has a crystal basis $(L,B(\la), l_{B(\la)})$
such that $L_{\la} / q L_{\la} = l_{b_{\la}}$.
Then we have

\bni

\item $\V \otimes M = \bigoplus_{\la + \epsilon_j : \text{strict}}
M_ j,$ where $M_j$ is a highest weight $\Uq$-module with highest
weight $\la + \epsilon_j$ and $\dim (M_{j})_{\la + \epsilon_j} = 2
\dim M_{\la}$,

\item $L_{j} = (\mathbf{L} \otimes L) \cap M_{j}$,

\item $\B \otimes B(\la) \simeq
\coprod_{\la+ \epsilon_j: \text{strict}} B_{j}$, where $$B_{j}
\simeq B(\la + \epsilon_j), \quad
L_j/qL_j =
\bigoplus_{b \in B_j} l_{b}.$$

\end{enumerate}
\end{enumerate}
 \enth

%\begin{proof}
We will prove all of our assertions at once by induction on
the length of $\la$. The proof is involved
because our theorem consists of several interlocking statements.
The key ingredient is a combinatorial proof of (c).
%\end{proof}

%\vskip 1cm

\end{document}